\documentclass[12pt]{amsart}
\usepackage{amssymb}
\usepackage{amsfonts}
\usepackage[cmyk]{xcolor}
\usepackage[margin=1.1in]{geometry}

\newcount\minleft
\newcount\timehour
\def\thetime{\timehour=\time
\divide\timehour by60 \minleft=\timehour \multiply\minleft by -60
\advance\minleft by\time \ifnum\time>720\advance\timehour
by-12\fi\relax
\number\timehour:\ifnum\minleft<10 %
    0\fi\relax\number\minleft
    \ifnum\time>720~pm \else~am\fi}

\newtheorem{theorem}{Theorem}[section]

\newtheorem{corollary}[theorem]{Corollary}
\newtheorem{proposition}[theorem]{Proposition}
\newtheorem{remark0}[theorem]{Remark}
\newtheorem{example0}[theorem]{Example}
\newtheorem{definition}[theorem]{Definition}

\newtheorem{question}[theorem]{Question}

\newenvironment{example}{\begin{example0}\rm}{\end{example0}}
\newenvironment{remark}{\begin{remark0}\rm}{\end{remark0}}

\newcommand{\propref}[1]{Proposition~\ref{#1}}
\newcommand{\thmref}[1]{Theorem~\ref{#1}}

\newcommand{\remref}[1]{Remark~\ref{#1}}

\def\max{{\mathfrak{m}}}                   

\def\maxn{{\mathfrak{n}}}                   
\def\res{{\mathbf{k}}}

\def\md{{{\bf 1}}}

\def\HF{{\operatorname{H\!F}}}

\def\deg{\operatorname{deg}}
\def\dim{\operatorname{dim}}
\def\ann{\operatorname{Ann}}

\def\length{\operatorname{Length}}

\def\socdeg{\operatorname{socdeg}}
\def\spec{\operatorname{Spec}}

\begin{document}
\title[Inverse systems of  one-dimensional Gorenstein   rings]{{\bf
A constructive approach to one-dimensional Gorenstein  $\res$-algebras}}
\author[J. Elias]{J. Elias ${}^{*}$}
\thanks{${}^{*}$
Partially supported by  MTM2016-78881-P}
\address{Joan Elias
\newline \indent Departament de Matem\`{a}tiques i Inform\`{a}tica
\newline \indent Universitat de Barcelona (UB)
\newline \indent Gran Via 585, 08007
Barcelona, Spain}  \email{{\tt elias@ub.edu}}

\author[M. E. Rossi]{M. E. Rossi ${}^{**}$}
\thanks{${}^{**}$
Partially supported by PRIN 2010-11
``Geometria delle varieta' algebriche".\\
\rm \indent 2010 MSC:  Primary
13H10; Secondary 13H15; 14C05}

\address{Maria  Evelina  Rossi
\newline \indent Dipartimento di Matematica
\newline \indent Universit\`{a} di genova
\newline \indent Via Dodecaneso 35, 16146 Genova, Italy}
\email{{\tt rossim@dima.unige.it}}

\date{\today}

\begin{abstract}
Let     $R$ be the power series ring or the  polynomial ring over a field $\res$     and let $I  $ be an ideal of $R.$ 
Macaulay proved that the Artinian Gorenstein $\res$-algebras $R/I$ are in one-to-one correspondence with the cyclic $R$-submodules  of the divided power series ring $\Gamma.  $  The result is effective in the sense that any polynomial of degree $s$ produces an Artinian Gorenstein $\res$-algebra of socle degree $s.$ In a recent paper,  the  authors   extended Macaulay's correspondence
 characterizing  the $R$-submodules of $\Gamma $  in one-to-one correspondence with Gorenstein d-dimensional $\res$-algebras. However, 
  these submodules
in positive dimension are not finitely generated. 
Our goal  is  to give constructive and finite procedures for the construction of Gorenstein $\res$-algebras of  dimension one and any codimension. This has been achieved through a deep analysis of the $G$-admissible submodules of $\Gamma. $   Applications to the Gorenstein linkage  of zero-dimensional schemes and to Gorenstein affine semigroup rings are discussed.
\end{abstract}

\maketitle

\section{Introduction }

Gorenstein rings were introduced by A. Grothendieck and they are a generalization of complete intersections,   indeed the two notions coincide in codimension two by a well known result by Serre.
Codimension three Gorenstein rings are completely described   by  Buchsbaum and Eisenbud's structure theorem, \cite{BE77}, but despite many attempts the  general construction of Gorenstein rings remains  an open problem in higher codimension.
A. Kustin and M. Miller, in a series of papers, studied the structure of  Gorenstein ideals of codimension $4$,  see \cite{KM83}. More recently, M. Reid  studied their
projective resolution   aiming to extend the result of D. Buchsbaum and D. Eisenbud, see  \cite{Rei13}. Gorenstein rings are of great interest in many areas of mathematics and they have appeared as an important component in a significant number of problems with  applications
to commutative algebra,  singularity theory, number theory  and more recently to combinatorics, among other areas. The lack of a general structure of   Gorenstein rings is the main obstacle in several problems, among them the Gorenstein linkage. For a complete and interesting presentation of Gorenstein rings, see H. Bass's   and  C. Huneke's papers,  \cite{Bas63}  and \cite{Hun99}.

Let  $R$ denote the power series ring $ \res[\![z_1, \dots, z_n]\!] $ or the  polynomial ring
$ \res[z_1, \dots, z_n] $ over a field $\res$ and let $I \subset R$ be an ideal (homogeneous if $R$ is a polynomial ring).
We denote by $\mathcal M$ the maximal  ideal of $R$ generated by $z_1,\cdots, z_n$.

As an effective consequence of Matlis duality,   it is known that an  Artinian ring $R/I $ is a Gorenstein $\res$-algebra if and only if $I$ is the set of solutions of a system of polynomial differential operators with constant coefficients.
Macaulay,  at the beginning of the 20th century,   proved that the Artinian Gorenstein $\res$-algebras are in correspondence with the cyclic $R$-submodules  of     the divided powers ring $\Gamma=\res_{DP}[Z_1,\cdots, Z_n] $ where the elements of $R$ act as derivatives (or contraction) on $\Gamma$, see \cite{Ems78}, \cite{IK99}.

Thanks to the effective construction of the Artinian Gorenstein algebras, in the last twenty years several authors have applied this device to several problems, among others:   Waring's problem  \cite{Ger96},
the n-factorial conjecture in combinatorics and geometry  \cite{Hai94},
the cactus rank \cite{RS13},  the geometry of the punctual Hilbert scheme of Gorenstein schemes \cite{IK99}, the classification up to analytic isomorphism of Artinian Gorenstein rings \cite{ER12}, and  the Koszulness of $\res$-algebras  \cite{CRV}.

Recently    in \cite{ER17} the authors extended   Macaulay's correspondence
 characterizing  the submodules of $\Gamma, $  called $G$-admissible submodules,  in one-to-one correspondence with Gorenstein d-dimensional $\res$-algebras
 (Theorem \ref{bijecGor}).   These submodules
in positive dimension are not finitely generated. Clearly this is  an obstacle to the effective construction of Gorenstein algebras of positive dimension.

In this paper we  give a constructive and finite procedure  for producing  Gorenstein one-dimensional $\res$-algebras of any codimension.    Thanks to a deep analysis of the structure of  the $G$-admissible $R$-submodules  of $\Gamma, $  it is possible to write an algorithmic procedure for constructing  step-by-step a finite subset of $\Gamma, $ called a $G$-admissible set, as a good candidate  for being extended to a $G$-admissible submodule.   In the graded case, similarly to the Artinian case,  a suitable  DP-polynomial  $H\in \Gamma$  uniquely determines a one-dimensional Gorenstein ring and hence a $G$-admissible $R$-submodule, see Proposition \ref{finite}, Remark \ref{homogeneous}.  This is no  longer true in the local case where more sophisticated techniques will be necessary.  
  By using the theory of the standard bases (or Hironaka bases),   Theorem \ref{lift} gives necessary and sufficient effective conditions to build up   the  Gorenstein one-dimensional $\res$-algebras  from  a finite subset of $\Gamma.$   In particular, starting from a suitable {\bf{finite}}  $G$-admissible set $\mathcal H$,  we can construct  all the   ideals $I $ of $R$ such that $R/I$ is a Gorenstein one-dimensional $\res$-algebra  and the dual  contains  $\mathcal H.$   Example \ref{non-iso} shows that two Gorenstein rings  sharing the same $\mathcal H$ in general are not  analytically isomorphic, even if they share the same associated graded ring.  In the last two sections we apply the results to  two classes of examples:  the construction of  the inverse system of  $d$ distinct Gorenstein points of ${\mathbb P^3_{\res}}$ and  the inverse system of a class of    Gorenstein semigroup rings
in ${\mathbb A^3_{\res}}$. The computations are performed  by using the computer algebra system  Singular  \cite{DGPS}  and the Singular library INVERSE-SYST.lib \cite{E-InvSyst14}.


\section{The structure of the Inverse system}

Let $V$ be a vector space of dimension $n$ over a field $\res$ where, unless specifically stated otherwise,  $\res$ is  an infinite  field of any characteristic.
If  $V $ denotes the $\res$-vector space $\langle z_1, \dots, z_n \rangle,$ then we denote  by   $V^* = \langle Z_1, \dots, Z_n \rangle$   the  dual basis.
Let $P= Sym_{\cdot}^{\res} V = \oplus_{i \ge 0} Sym_i^{\res }V $ be  the  standard graded polynomial ring in $n$ variables over $\res$ and $\Gamma = D_{\bf \cdot}^{\res}(V^*) = \oplus_{i\ge 0} D_i^{\res}(V^*) = \oplus_{i\ge 0}  {\rm Hom}_{\res} (P_i, \res) $    be the graded $P$-module of graded $\res$-linear homomorphisms from $P$ to $\res, $ hence
   $\Gamma \simeq  \res_{DP} [Z_1, \dots, Z_n] $ the divided power ring.
   In particular   $\Gamma_j = \langle \{Z^{[L]} \  | \  |L|=j \} \rangle $ is the span of the dual generators to $z^{L}=z_1^{l_1}\cdots z^{l_n}$ where $L$ denotes the multi-index $L=(l_1,\dots, l_n)\in \mathbb N^n $ of length $|L|= \sum_i  l_i.$
If $L\in \mathbb Z^n$ then we set $ Z^{[L]}=0$ if any component of $L$ is negative.   The monomials  $ Z^{[L]} $ are called divided power monomials (DP-monomials) and the elements $F= \sum_L b_L Z^{[L]} $ of $\Gamma $ the divided power polynomials (DP-polynomials).

We recall   that $\Gamma $ is an  $R$-module with $R$ acting   on  $\Gamma$ by {\it {contraction}}   as follows. This action is also called apolarity. {

\begin{definition} If $h =\sum_M a_M z^M \in R $ and $F= \sum_L b_L Z^{[L]} \in \Gamma, $ then the contraction of $F$ by $h$ is defined as
$$ h \circ F =  \sum_{M, L} a_M b_L Z^{[L-M]} $$
\end{definition}

\noindent
   For short,   from now on we write $Z^L$ instead of $Z^{[L]}$.

Recall that the injective hull $E_R(\res) $ of $\res$ as an $R$-module is isomorphic as an $R$-module to the divided power ring $\Gamma$ (see \cite{Gab59}, \cite{Nor72b}).
For detailed information   see  \cite{EisGTM}, \cite{Ems78}, \cite{IK99}, Appendix A.
If the characteristic of the field $\res$ is zero, then there is a natural isomorphism of $R$-algebras between $(\Gamma , \circ) $,
where $\circ$ is the contraction already defined, and the usual polynomial ring $P $  replacing   contraction with taking     partial derivatives.      In this paper    we  always consider  $\Gamma$ as an $R$-module by  contraction.

If $I\subset R$ is an ideal of $R$ then  $(R/I)^\vee={\rm Hom}_R(R/I,E_R(\res))$ is the $R$-submodule of $\Gamma$
$$
{I^{\perp}} =\{ g \in \Gamma \ |\  I \circ g = 0 \  \}.
$$
This submodule of $\Gamma $  is called {\it{Macaulay's inverse system of $I$.}}
Given an $R$-submodule $W$ of $\Gamma, $ then the dual $W^\vee={\rm Hom}_R(W,E_R(\res))$ is the ring $R/\ann_R(W)$ where
$$ \ann_R(W)= \{ f \in R \ \mid \ f \circ g= 0 \ \mbox{ for \ all \ } g \in W\} 
$$
 is an ideal of $R$.
If $I$ is a homogeneous ideal of $R$ (resp. $W$ is generated by homogeneous polynomials) then $I^{\perp}$ is generated by homogeneous polynomials of $\Gamma$ (resp. $\ann_R(W) $ is an homogeneous ideal of $R$)  and
$ I^{\perp} = \oplus I_j ^{\perp} $ where $I_j ^{\perp} = \{ F \in \Gamma_j \mid  \ h \circ F=0 \ \ \mbox{ for \ all \ }  h \in I_j \}.$ Notice that the Hilbert function of $R/I$ can be computed by $I^{\perp}, $ see for instance \cite{ER17} Section 2.
 Macaulay in \cite[IV]{mac16} proved a particular case of Matlis duality, called Macaulay's correspondence,
between the
  ideals $I\subseteq R  $ such that $R/I$ is an Artinian local ring and $R$-submodules  $W=I^{\perp}$ of $\Gamma $ which are finitely generated.
Macaulay's  correspondence is an effective method for computing Gorenstein Artinian rings, see  \cite{CI12}, Section 1, \cite{Iar94}, \cite{Ger96} and \cite{IK99}. We summarize with the statement 
 {\it{Artinian  Gorenstein $\res$-algebras $A=R/I$ of socle degree $s$ correspond to cyclic $R$-submodules  of $\Gamma$ generated by a DP-polynomial $F\neq 0$ of degree $s$.}}


The authors extended  Macaulay's correspondence characterizing the $d$-dimensional local Gorenstein $\res$-algebras in terms of suitable submodules of $\Gamma$, see \cite{ER17}. This result in the one-dimensional case will be our starting point and for completeness we include here the statement.

\begin{definition}
\label{G-modules}  An   $R$-submodule $M $ of  $\  \Gamma  $    is called $G$-admissible if it admits a countable system of  generators $\{H_l\}_{ l\in \mathbb N_+} $ satisfying  the following conditions
\begin{enumerate}
\item[(1)] There exists a linear form   $z \in R$ such that for all $l\in \mathbb N_+$
$$
z \circ H_l =
\left\{
\begin{array}{ll}
 H_{l-1}   & \text{ if }  l> 1 \ \\
0 & \text{ otherwise.}
\end{array}
\right.
$$
\item[(2)] $\ann_R(\langle H_l\rangle)\circ H_{l+1}=\langle H_{1}\rangle$ for all $l\in \mathbb N_+$.
\end{enumerate}

\noindent
If this is the case, we say that $M = \langle H_l, l\in \mathbb N_+\rangle$ is a
$G$-admissible $R$-submodule of $\Gamma$ with respect to the linear form   $z \in R$.
\end{definition}

Notice that if $R/I$ has positive depth, then  there always exists a linear form $z \in R$ which is regular modulo $I $ because $\res$ is infinite.

\begin{theorem}
\label{bijecGor}
There is a one-to-one correspondence $\mathcal C$ between the following sets:

\noindent
(i)
 one-dimensional  Gorenstein $\res$-algebras $A=R/I$,

\noindent
(ii) non-zero $G$-admissible $R$-submodules $M=\langle H_l, l\in \mathbb N_+\rangle$ of $\Gamma$.

\noindent
In particular, given an ideal $I\subset R$ with $A= R/I$ satisfying $(i)$ and $z$ a linear regular element  modulo $I, $   then
 $$\mathcal C(A)= I^{\perp} = \langle H_l,  l\in \mathbb N_+\rangle \subset S   \ \ {\text{with}}\ \
  \langle H_l\rangle=(I+(z^l))^\perp $$ is $G$-admissible.
 Conversely, given an  $R$-submodule $M$ of $\Gamma $ satisfying  (ii), then $$\mathcal C^{-1}(M)=R/I  \ \ {\text{with}}\ \
 I= \ann_R(M) = \bigcap_{l\in \mathbb N_+}\ann_R(\langle H_l\rangle).$$
\end{theorem}

The main goal of this paper is to construct     $G$-admissible $R$-submodules of $\Gamma.$
\medskip

\noindent Notice that  if we fix any DP-polynomial $H_1\in \Gamma  $ such that $z \circ H_1=0, $ then $M=\langle Z^l H_1, l\in \mathbb N \rangle$  is  $G$-admissible with respect to $z,$  but the corresponding one-dimensional Gorenstein ring is not of great interest because it is a cone over  the Gorenstein Artinian ring corresponding to the DP-polynomial $H_1, $   see \cite{ER17}, Proposition 4.1. In the next remark  we observe that  the choice of $H_1 $ is very important  in the construction of a $G$-admissible module, because it encodes much  information on the corresponding Gorenstein ring.

\begin{remark}
\label{details}
Let $A=R/I $ be a $1$-dimensional Gorenstein  ring and let $z$ be a linear regular element  modulo $I. $
Let  $I^{\perp} =\langle H_l, l\in \mathbb N_+\rangle$ be   the corresponding $G$-admissible dual module   with respect to $z.$   We recall that by Theorem \ref{bijecGor}, we have  $\langle H_1\rangle = (I+(z))^{\perp}. $
Hence by Macaulay's correspondence,  $\deg(H_{1})$ coincides with the  socle degree of the Artinian reduction $A/ z A.$
On the other hand we have the following inequality on the multiplicity $e(A)$ of $A: $
$$
e(A)\le \length_R(A/z A)=\dim_{\res}(\langle H_{1}\rangle).
$$
In particular the  equality holds  if $z$ is a superficial element  of $A. $
If  $A$ is a standard graded $\res$-algebra which is Gorenstein, then the above equality holds   for every regular element $z$ modulo $I. $    Hence $e(A) = \dim_{\res}(\langle H_{1}\rangle) $ (the dimension as a $\res$-vector space of the $R$-module generated by $H_1$) and $reg(A) = \deg H_{1} $ where $reg(A)$ is  the Castelnuovo-Mumford regularity of $A. $  

\end{remark}

As a consequence, in the graded case,  important geometric information such as the multiplicity, the arithmetic  genus, or more  generally,  the Hilbert polynomial of the Gorenstein $\res$-algebra can be controlled by the choice of $H_{1}, $ the first step in the construction of a  $G$-admissible dual module.

Recall that if  $I$ is an homogeneous ideal, then the  dual $R$-submodule $I^{\perp}$    of $\Gamma $ can  be generated by homogeneous (in the usual meaning) DP-polynomials and conversely if the $R$-submodule of $\Gamma$ is homogeneous, then the ideal $I$ is homogeneous.
Hence the correspondence will be between one-dimensional Gorenstein standard graded $\res$-algebras  and $G$-admissible homogeneous $R$-submodules    of $\Gamma.  $ The following result  refines Theorem \ref{bijecGor} in the case of graded $\res$-algebras.

\begin{theorem}
\label{gradedGor}
There is a one-to-one correspondence $\mathcal C$ between the following sets:

\noindent
(i) one-dimensional {\bf{Gorenstein standard  graded $\res$-algebras}}  $A=R/I$  of multiplicity $e=e(A) $ (resp. Castelnuovo-Mumford regularity $r= reg(A)$)

\noindent
(ii) non-zero {\bf{$G$-admissible homogeneous $R$-submodules}} $M=\langle H_l, l\in \mathbb N_+\rangle$ of $\Gamma$  such that $\dim_{\res} \langle H_{1}\rangle = e$ (resp. $\deg H_{1} =r$)
\end{theorem}

In the construction of $G$-admissible $R$-modules,  the following remark suggests to us to proceed step by step starting from $H_1.$

\begin{remark}
\label{integral}   Given two  DP-polynomials $H, G\in \Gamma, $ we say that $G$ is a primitive of $H$ with respect to $z  \in R  $ if
$z \circ G=H$.
From the definition of the contraction $\circ, $ we will get
$$G= Z   H + C$$
for some $C\in \Gamma$ such that $z  \circ C=0.  $   We remark that $Z   H $  denotes the usual multiplication in a polynomial ring and we do not use the internal multiplication in $\Gamma$ as DP-polynomials.
Hence, according to Definition \ref{G-modules}, part (1),  given  a $G$-admissible module generated by $\{H_l\}_{ l\in \mathbb N_+} $, we have that $H_{l+1}$    is a primitive  of $H_{l}$
with respect to $z   $ for every positive integer $l.  $ In particular $H_2$ is a primitive of $H_1$ and so on. Hence there exist $C_1, \dots, C_l $ in $\Gamma $ (depending on condition (2) of Definition \ref{G-modules}) such that $z \circ C_i =0$ and   for all positive $l$
\begin{equation} \label{C}  H_{l+1} = z^l H_1 + z^{l-1} C_1 + \dots + z C_{l-1} + C_l  \end{equation}
In the effective construction of a Gorenstein ring $R/I, $ note that we may assume $\Gamma=\res_{DP}[Z_1, \dots, Z_n] $  and take $z=z_1.$
This presentation has the advantage that we may assume $H_{1}, C_1, \dots, C_l$  in $\res_{DP}[Z_{2}, \dots, Z_n].$
\end{remark}

\medskip
According to Remark \ref{integral}, starting from a  DP-polynomial $H_1,$  in each step we have to choose  the ``constants" $C_i  $ imposing condition (2) of Definition \ref{G-modules}. We have implemented a Singular routine  to determine the possible DP-polynomials $H_{l+1} $ from $H_l$, if they exist \cite{DGPS}.


Guided  by the construction of  one-dimensional Gorenstein $\res$-algebras $R/I$  in a finite and effective number of steps,  we state  the following questions:

\begin{question} \label{2} When can a  finite set  $\{H_1, \dots, H_{t+1}\} $ of DP-polynomials of $\Gamma$  verifying   conditions (1) and (2) of Definition \ref{G-modules}  be completed to a system of generators of a $G$-admissible $R$-module $M=\langle H_l, l\in \mathbb N_+\rangle$?
\end{question}

\begin{question} \label{1}Does there exist an integer $t$ such that the finite subset $\{H_1, \dots, H_{t+1}\} $ of a $G$-admissible $R$-module $M=\langle H_l, l\in \mathbb N_+\rangle$  determines  (uniquely?) the Gorenstein ideal $I= Ann_R(M)$?
\end{question}

\medskip
According to the definition of a $G$-admissible $R$-submodule of $\Gamma,$  see Definition \ref{G-modules}, we say:

\begin{definition} A  {\bf finite set}  of DP-polynomials  $\mathcal H=\{H_l \ | \  0 < l \le t_0\}$  is a {\bf G-admissible set}  with respect to a linear element $z \in R$ if it satisfies  conditions (1) and (2) of Definition \ref{G-modules}  with respect to $z \in R. $
\end{definition}

Note that an algorithm in Singular has been implemented for computing $G$-admissible sets with respect to $z \in R$ starting from a $DP$-polynomial $H_1$ such that $z \circ H_1=0.$

\medskip

\begin{remark}
\label{(1)vs(2)}
Notice that in $\mathcal H=\{H_l; 0 < l \le t_0\}$ the last polynomial $H_{t_0}$ determines the full sequence $H_1,\cdots , H_{t_0-1}$ by contraction with respect to $z$. Thus  $\mathcal H$ is identified by the DP-polynomial $H_{t_0}. $
But $H_{t_0}$ is not an arbitrary polynomial simply of the shape described in Remark \ref{integral}.
In fact, in Definition \ref{G-modules}, condition $(1)$  does not imply condition $(2)$ as the following example shows.
Let us consider the finite set:
$\mathcal H=\{H_1=X^2, H_2=X^2Y, H_3=X^4+X^2Y^2, H_4=X^4Y+X^2Y^3, H_5=X^4Y^2+X^2Y^4\}$.
This set satisfies condition $(1)$ with respect to $y$;  and it satisfies condition $(2)$ except for the last polynomial since
$\ann_R(\langle H_4\rangle)\circ H_{5}   \supsetneqq \langle H_{1}\rangle$.
However, if we replace $H_5$ with  $H_5+X^6$ we get a G-admissible set with respect to $y.$  \end{remark}

From now on if $J$ is an ideal of $R$ (not necessarily homogeneous), we denote by $J_{\le s} R$ the ideal of $R$ generated by all the elements of $J$ in $\mathcal M^s \setminus \mathcal M^{s+1} $ if any, otherwise $(0).$ In the homogeneous case, it means generated by forms of $J$ of degree at most $s.$  The next result gives a first positive answer to Question \ref{1} in the homogeneous case.

\begin{proposition} \label{finite}
Let  $t_0$ be a positive integer and   let $\mathcal H=\{H_l \ |\  0 < l \le t_0\}$  be  a   homogeneous G-admissible  set with respect to $z \in R. $  Let   $ \deg H_1=r    $ and assume  $t_0 \ge r+2.$  If $\mathcal H $ can be extended to  a $G$-admissible $R$-module, then the corresponding graded Gorenstein   $\res$-algebra  $A=R/ I$ is uniquely determined and  $$ I= \ann_R(H_{r+2 })_{\le r+1}R.$$
\end{proposition}

 The previous result was proved in \cite{ER17}, Proposition 4.2.  It depends on  the   fact that   the maximum degree of a minimal system of generators of $I$ is at most  $ reg(A) +1$   and it  coincides with the socle degree $\socdeg(A/zA) +1= r+1.$ Clearly it is also true that  $$ I= \ann_R(H_{t_0 })_{\le r+1}R.$$
  The previous result can be improved if we know the maximum degree $t$ of the generators of $I $ in which  case one can replace   $r = \deg H_{\md} $ by $  t-1 (\le r),  $ hence
  $$ I= \ann_R(H_{t+1})_{\le t}R.$$

\begin{remark} As in the Artinian case, the previous result tells us that there is a one-to-one correspondence between Gorenstein graded $\res$-algebras of dimension one  and  socle degree $r$ and suitable homogeneous DP-polynomials $H_{r+2} $ of degree $2r +1.$
\end{remark}

We present here an example.

\begin{example} \label{no}
We consider the following $G$-admissible set of DP-polynomials in $\Gamma= \res_{PD}[X,Y]   $ with respect to $y$ of  $R=\res[\![x,y]\!]$:
$$
\mathcal H=\{ H_1=X^2, H_2=X^2Y, H_3=X^4+X^2Y^2, H_4=X^4Y+X^2Y^3\}
$$
It is a subset of the $G$-admissible $R$-submodule $M$ of $\Gamma$ generated by $$\{H_t= \sum_{i=0}^{2 i \le t-1}X^{2+2i}Y^{t-1-2i} \ |\ t\ge 1\}$$
\noindent In this case $r=\deg H_1=2, $   hence according  to Proposition \ref{finite}
$$ I= \ann_R(H_{4 })_{\le 3}R.$$
In fact, $I=\ann_R(H_{4 })_{\le 3}R =(x^3-x y^2) $   is  a Gorenstein ideal, $y$ is a non-zero divisor in $A=R/I$ and
 $I^{\perp} = M.$
 We remark that in particular $\mathcal H$ (actually $H_4$) determines the Gorenstein ideal $I$ and hence the $G$-admissible $R$-module $I^{\perp}.$
  \end{example}

We end this section with a very partial  answer (for the moment) to Question \ref{2} in the homogeneous case.  We will need the following proposition. 

\begin{proposition}
\label{tata}
Let $H=\{H_l; 0 < l \le t_0\}$ be a  {\it {finite}}   G-admissible set with respect to $z \in R$. Then
$$
\ann_R(H_{t+1})+(z^t)=\ann_R(H_t)
$$
for all $t=1,\cdots, t_0-1$.
In particular
$\ann_R(H_{t+1})+(z)=\ann_R(H_1).$
\end{proposition}
\begin{proof}
Since $z\circ H_{t+1}=H_t$ we have $\langle H_t\rangle \subset \langle H_{t+1}\rangle$, and then
$\ann_R(H_{t+1})\subset \ann_R(H_t)$.
From condition $(1)$ we deduce that $z^t\circ H_t=0$ so $(z^t)\subset  \ann_R(H_t)$.

Let now $a\in  \ann_R(H_t)$, from  conditions $(2)$ and $(1)$ there exists $b\in R$ such that
$$
a\circ H_{t+1}=b\circ H_1= b\circ (z^t\circ H_{t+1}).
$$
Hence $a- b z^t\in  \ann_R(H_{t+1})$, so $a\in  \ann_R(H_{t+1})+(z^t)$.
\end{proof}
\begin{corollary} \label{finite2} Assume $\dim R=2. $   Let  $t_0$ be a positive integer and let $\mathcal H=\{H_l; 0 < l \le t_0\}$ be a  G-admissible homogeneous set with respect to $z \in R. $
Assume  $t_0\ge r+2 $ where  $r = \deg H_{1} $.
Then  $\mathcal H $ can be extended to  a $G$-admissible module  and the corresponding graded Gorenstein   $\res$-algebra  $A=R/ I$ is uniquely determined by   $$ I= \ann_R(H_{r+2 })_{\le r+1}R.$$
\end{corollary}
\begin{proof}  Denote by $(x,z)$ a minimal system of generators  of the maximal ideal of $R. $ Because $z \circ H_1=0 $ and $r = \deg H_{1}, $ we may write   $H_1=x^r$ and hence $\ann(H_1)=(x^{r+1},z).$  By Proposition \ref{tata} we have $\ann(H_{t+1}) +(z)= \ann(H_1) $ for every $t\ge 1,$  hence $\ann(H_{t+1}) \subseteq (x^{r+1}, z).$  We also know that $\ann(H_{t+1}) $ is generated by a regular sequence, say  $(F_t,  G_t). $ Hence we may assume $F_t= x^{r+1} +zM_t$ for some form $M_t$ of degree $r $ and $G_t =z N_t $ for some form $N_t \in R.$ We prove $\deg G_t \ge t+1 $ for every $t=0,\dots, t_0. $   We proceed inductively on $t.$ If $t=0, $ then $M_0=0 $ and $G_0=z. $ Assume $ \ann(H_{j+1}) =(x^{r+1}+zM_j, z N_{j}) $  with $\deg N_j \ge j$ and prove that   $ \ann(H_{j+2}) =(x^{r+1}+ zM_{j+1}, zN_{j+1}) $  with $ \deg N_{j+1} \ge j+1. $ Now  $z N_{j+1}  \circ H_{j+2} =0   $ implies   $N_{j+1}  \circ H_{j+1}=0.   $ Hence  $N_{j+1}  \in (x^{r+1}+zM_j, z N_{j}) $ and  we conclude. Hence $\ann_R(H_{r+2 })_{\le r+1}R= ( x^{r+1}+ zM_{r+1}) = I. $
\end{proof}

Unfortunately,   \propref{finite} and \propref{finite2}  cannot be extended to  the non-homogeneous case, even in the algebraic case.   In fact the information on the multiplicity is not enough to determine the singularity.  The following example is a counterexample in the non-graded case to both the above results.

\begin{example}
For all $n\ge 2$ we consider the one-dimensional local ring $A_n=k[\![x,y]\!]/(f_n)$ with
$f_n=y^2-x^n. $ Notice that, for  all $n\ge 1,$  $A_n$ is an algebraic  Gorenstein ring of multiplicity $e(A_n)=2.$
On the other hand
$A_n/x A_n=k[y]/(y^2) $ is an Artinian reduction of $A_n$, so $H_1=Y$ and hence $\deg H_1=1.$
If $n \ge 4, $  it is clear that we cannot obtain  the ideal $(f_n)$ after $\deg H_1 +2=3 $ steps as \propref{finite} and Corollary \ref{finite2}  suggest.
\end{example}

\section{One dimensional Gorenstein local or graded rings}

The aim of this section is to give constructive answers to Questions \ref{1} and \ref{2} in the  local and graded case.  It  will be useful to recall some well known facts concerning the standard bases in the local case. Let  $R $ be   the ring of  formal power series with coefficients in $\res$   with  maximal ideal $\mathcal M   $ and let $P $ be  the corresponding polynomial ring. Notice that $P=gr_{\mathcal M }(R) $ is the associated graded ring of $R$.
 For every element   $f\in R\setminus \{0\}$ we can write $f=f_v+f_{v+1}+\cdots$, where $f_v$ is not zero and $f_j$ is a homogeneous polynomial of degree $j$ in $P$ for every $j\ge v.$
 We say that $v(f) $ is the order of $f$ (or $\mathcal M$-valuation), denote $f_v$ by $f^*$ and call it the initial form of $f.$
 If $f=0$ we agree that its order is $\infty.$ We denote by $I^*$  the ideal of $P$ generated by all the initial forms of the elements in $I.$ A  set of generators (not necessarily minimal) of $I, $ say $\{f_1, \dots, f_r\}, $ is called   a  {\bf{standard basis}} (or Hironaka's basis)   of $I$ if their initial forms generate $I^*.$ Notice that the associated graded ring of $A=R/I$ with respect to the maximal ideal is isomorphic to $P/I^*.$  It is known  that the Krull dimension of $A$ and of $gr_{\mathcal M }(R/I) =P/I^*$ coincide, as well, by definition,  their Hilbert function. Unfortunately, even if $A$ is Gorenstein, then in general  $gr_{\mathcal M }(R/I) $ is no longer Gorenstein nor even  Cohen-Macaulay, see for instance the examples in Section 5. Notice that if $I$ is a homogeneous ideal in $P,$ then   any homogeneous system of generators of $I$ is a standard basis.
A  set of generators  of $I $   is called    a {\bf{minimal standard basis}}  of $I  $ if the  initial forms minimally generate   $I^*.$  Notice that the orders  of the elements in a minimal system of generators of $I$ are not uniquely determined, instead the orders of the elements of a minimal standard basis of $I$ (and the number of generators)   are uniquely  determined by the first graded Betti numbers of $ I^*.$     If $\{f_1, \dots, f_r\}  $ is   a   minimal  standard basis   of $I, $ then $v(f_i) = \deg(f_i^*) $ and $r$ is the minimal number of generators of the homogeneous ideal  $I^*.$ We recall that an element $z \in  \mathcal M\setminus \mathcal M^2 $ is superficial for $A$ if and only if $z^*$ is a homogeneous filter-regular element of degree one in $G= gr_{\mathcal M }(R/I) =P/I^*$ or equivalently $(0:_G z^*)_j=0$ for all large degree $j.$ If  depth $A>0,$ then a superficial element is also a regular element modulo $I,  $ see for instance \cite{RV-Book}, Section 1.2.

\vskip 2mm
The next Theorem has a central role in this paper.  In particular it extends and it refines Proposition \ref{finite} in the local case.

\begin{theorem}
\label{minimalGB}
Let $A=R/I$ be a one-dimensional Gorenstein ring  and let
$I^{\perp}=  \langle  H_i \mid i \in \mathbb N_+ \rangle$ be the correspondent   $G$-admissible module with respect to $z$  (a regular linear element  modulo  $I$).  Given $t\ge e= \dim_{\res}\langle H_1 \rangle  $, let $\{h_1,\dots, h_r\} $ be the elements of a minimal   standard basis  of $\ann_R(H_{t+1}) $ such that $v(h_i) \le t$.
 Then
\begin{enumerate}
\item[(i)]  $(\ann_R(H_{t+1})^*)_{\le t} = (I^*)_{\le t}$.
\item[(ii)] There exist $\alpha_1, \dots, \alpha_r \in R$
  such that  $\{h_1+ \alpha_1 z^{t+1}, \dots, h_r+ \alpha_r z^{t+1}\}$
is a minimal  standard    basis of $I$.
\item[(iii)] $\ann_R(H_{t+1})=I+(z^{t+1})=(h_1,\dots,h_r)+(z^{t+1})$.
\end{enumerate}
Assume that $z$ is a superficial element    of $A$ and  denote by $J$ the ideal generated by $\{h_1,\dots, h_r\} $ and $B=R/J$.
\begin{enumerate}
\item[(iv)]
 $B$ is a one-dimensional Gorenstein ring and  $z$ is a regular element  modulo $J.$ In particular
$J^*=I^*$ and hence the Hilbert functions of $A$ and $B$ agree.
\item[(v)]
Assume that $\{h_1,\dots,h_s\}$ is a minimal system of generators of  $J$,
then $\{h_1+ \alpha_1 z^{t+1}, \dots, h_s+ \alpha_s z^{t+1}\}$
is a minimal system of generators of $I$.
\end{enumerate}
\end{theorem}
\begin{proof}
$(i)$
Since $I+ (z^{t+1})=\ann_R(H_{t+1})$ we get that $I^* \subset \ann_R(H_{t+1})^*$. Conversely,
let $F\in \ann_R(H_{t+1})^*_{\le t}$ be a homogeneous form.
Then there are $f\in I$ and $\beta\in R$ such that $(f+\beta z^{t+1})^*= F$.
Since the degree of $F$ is less than or equal to $t, $ we deduce
that $F=f^*\in I^*_{\le t}$.

\noindent
$(ii)$
Let $g_1,\dots, g_p$ be a minimal   standard basis of $I$.
Since $t\ge e= \dim_{\res}\langle H_1 \rangle \ge e(R/I)$ and  $I^*$ is minimally generated by forms of degree less than  or equal to $e(R/I)$,
\cite{Eli86c}, we have that $\deg g_i^*\le t$, $i=1,\cdots, p$.

Let $\{h_1,\dots, h_r\} $ be the elements of a minimal  standard basis of $\ann_R(H_{t+1}) $ such that $\deg h_i^* \le t.$
From $(i)$ we get
$$
(g_1^*,\dots, g_p^*)_{\le t}=
I^*_{\le t}=
\ann_R(H_{t+1})^*_{\le t}=
(h_1^*,\dots, h_r^*)_{\le t}
$$
and then
\begin{equation}
\label{minstd}
I^*=
(g_1^*,\dots, g_p^*)=
(h_1^*,\dots, h_r^*).
\end{equation}
Thus  $h_1^*,\dots, h_r^*$ is a minimal system of generators of $I^*$.

For all  $h_i\in I+ (z^{t+1})$ we can write, $i=1,\dots,r$,
$$
h_i=\gamma_i -\alpha_i z^{t+1}
$$
with $\gamma_i\in I$ and $ \alpha_i\in R$.
From this identity we have $h_i+ \alpha_i z^{t+1}=\gamma_i \in I$, so
$$
h_i^*=(h_i+ \alpha_i z^{t+1})^*\in I^*
$$
for $i=1,\dots,r$.
Hence  $\{h_i+ \alpha_i z^{t+1}\}_{i=1,\dots,r}$ is a minimal   standard   basis of $I$.

\noindent
$(iii)$
It is a consequence of $(ii)$.

\noindent
$(iv)$
We consider the morphism $\pi$
$$
\pi: X=\spec\left(\frac{R[w]}{(h_1+ w^{t+1} \alpha_1 z^{t+1}, \dots, h_r+ w^{t+1} \alpha_r z^{t+1})}\right)  \longrightarrow \spec(\res[w])
$$
Notice that $\pi^{-1}(0)$ is $B$ and that for all $a\neq 0\in \res$ we have
$\pi^{-1}(a)\cong A$.

First we prove that $B$ is one-dimensional.

We assume that $z$ is a superficial element of degree one of $A$,
so $e$ is the multiplicity of $A$.
We denote by $\max$ (resp. $\maxn$) the maximal ideal of $A$ (resp. $B$).

From the identity
$I+(z^{t+1})=J +(z^{t+1})$, $t\ge e$,
we deduce $I+\mathcal M^{e+1}=J+\mathcal M^{e+1}$,   so  $\max^{j}/\max^{j+1}= \maxn^{j}/\maxn^{j+1}$ for $j=e-1,e$.
Since the multiplication  by $z$ induces an isomorphism
$$
\frac{\max^{e-1}}{\max^{e}} \overset{.z}{\longrightarrow} \frac{\max^{e}}{\max^{e+1}}
$$
we get an isomorphism
$$
\frac{\maxn^{e-1}}{\maxn^{e}} \overset{.z}{\longrightarrow} \frac{\maxn^{e}}{\maxn^{e+1}}
$$
By Nakayama's lemma we deduce $\maxn^{e}= z \maxn^{e-1}$, so
$\maxn^{j}= z \maxn^{j-1}$ for all $j\ge e$.
\vskip 2mm
From this we get that $\dim(B)\le 1$.
Assume that $\dim(B)=0$. From the upper semi-continuity of the dimension over the ground field $\res$ of the fibers of the morphism $\pi$ we get $\dim(A)=0$ which is  a contradiction with the hypothesis
$\dim(A)=1$.

From the identity
$I+(z^{t+1})=J+(z^{t+1})$ we get
$I+(z)=J+(z)$, so the multiplicity $e'$ of $B$ satisfies
$$
e'\le \dim_{\res}(R/J+(z))=\dim_{\res}(R/I+(z))=e.
$$

We denote by $e' T -e'_1$ (resp. $e T -e_1$) the  Hilbert Polynomial  of $B$ (resp. $A$).
Assume  $e'< e$ and let $l$ be an integer such that
$$
e' l-e'_1< e l- e_1.
$$
Consider now the morphism
$$
\pi_l: X_l=\spec\left(\frac{R[w]}{(h_1+ w^{t+1} \alpha_1 z^{t+1}, \dots, h_r+ w^{t+1} \alpha_r z^{t+1})+\mathcal M^{l+1}R[w]}\right)  \longrightarrow \spec(\res[w])
$$
We have $\dim_{\res}(\pi_l^{-1}(0))=e' l-e'_1$ and that
$\dim_{\res}(\pi_l^{-1}(a))=e l- e_1$ for all $a\neq 0\in \res$.
By upper semi-continuity of the dimension over the ground field $\res$ of the fibers of the morphism $\pi_l$ we get a contradiction, so $e'=e$.

We know that $\maxn^{j}= z \maxn^{j-1}$ for all $j\ge e=e'$
so
$$
\frac{\maxn^{j-1}}{\maxn^{j}} \overset{.z}{\longrightarrow} \frac{\maxn^{j}}{\maxn^{j+1}}
$$
is an isomorphism for all $j\ge e'=e$.
Hence $z$ is a superficial degree one element of $B$.
Since $\dim_{\res}(B/(z))=e'$ we get that $z$ is a nonzero divisor of $B$ by \cite{RV-Book}, Proposition 1.2(3).

Since $A$ and  $B$ are  one-dimensional Cohen-Macaulay local rings of multiplicity $e$,
$I^*$ and $J^*$ are minimally generated by forms of degree less than or equal to $e$.
On the other hand
$$
J^*_{\le e}=(J+(z^{t+1}))^*_{\le e}=(I+(z^{t+1}))^*_{\le e}=(h_1^*,\cdots ,h_r^*)_{\le e}=I^*_{\le e}
$$
Hence $I^*=J^*$ and the Hilbert functions of $A$ and $B$ agree.

\noindent
$(v)$
Assume that $\{h_1,\cdots, h_s\}$, $s\le r$, is a minimal system of generators of the ideal
$J$.
Assume now that $s< r$ and let $s< j\le r$.
Then there exist $a_{j}^i\in R, i=1,\cdots, s$,
such that
$$
h_j=\sum_{i=1}^{s} a_j^i h_i
$$
Thus 
$$
h_j+\alpha_j z^{t+1}- \sum_{i=1}^{s} a_j^i (h_i+ \alpha_i z^{t+1})=
z^{t+1}(\alpha_j- \sum_{i=1}^{s} a_j^i\alpha_i)\in I
$$
Since $z$ is a non-zero divisor modulo $I$ we deduce
$$
\alpha_j- \sum_{i=1}^{s} a_j^i\alpha_i=\beta_j\in I
$$
and then
$$
h_j+\alpha_j z^{t+1}=\sum_{i=1}^{s} a_j^i (h_i+ \alpha_i z^{t+1})+ \beta_j z^{t+1}
$$
for $j=s+1,\cdots, r$.
Since $\beta_j z^{t+1}\in \mathcal M I$, by  Nakayama's Lemma $I$ is generated  by
$\{h_i+\alpha_i z^{t+1}\}_{i=1,\cdots ,s}$.

Assume that $\{h_i+\alpha_i z^{t+1}\}_{i=1,\cdots ,s}$ is not a minimal system of generators,
for instance if
$$
h_s+\alpha_s z^{t+1}=\sum_{i=1}^{s-1} a_j^i (h_i+ \alpha_i z^{t+1})
$$
for some $a_{j}^i\in R, i=1,\cdots, s$, then
$$
h_s -\sum_{i=1}^{s-1} a_j^i h_i=   z^{t+1}(- \alpha_s+\sum_{i=1}^{s-1} a_j^i \alpha_i)
$$
Since  $z$ is a non-zero divisor modulo the ideal $J$  we get
$$
 \alpha_s-\sum_{i=1}^{s-1} a_j^i \alpha_i=\beta_j\in J
$$
and then
$$
h_s - \sum_{i=1}^{s-1} a_j^i h_i\in \max J
$$
which is a contradiction with the assumption that  $\{h_1,\cdots, h_s\}$ is a minimal system of generators of $J$.
\end{proof}

Notice that in  Theorem \ref{minimalGB}, the ideal  $\ann_R(H_{t+1}) = I +(z^{t+1})$  has generators  of valuation    $ v \le t  $ because  $ t \ge e(R/I).  $
Hence the assumptions of this result are consistent.

\begin{remark}
\label{homogeneous}   If the ideal $I$ is homogeneous, then $\ann_R(H_{t+1}) = \ann_R(H_{t+1})^* $ and $I=I^*. $ In this case, in the previous result,  we have  necessarily $\alpha_i=0$ for every $i.$  Hence $H_{e+1} $ determines the ideal  $I $ where $e=e(R/I).$   This confirms Proposition \ref{finite} because $ e\ge r+1 $ where $r =\deg H_1. $  We remark that  in the homogeneous case,  $H_{r+2} $ determines the Gorenstein ideal $I, $ while instead  in the non homogeneous case,  the DP-polynomial   $H_{e+1} $ determines a Gorenstein  ideal $I $  (but not   uniquely).   We will see that  $H_1, \dots, H_{e+1} $ can be completed in different ways to a $G$-admissible module, see Example \ref{2ideals}.
\end{remark}

The next result gives a sufficient condition for obtaining a Gorenstein 1-dimensional ring from a finite sequence $\mathcal H $ only requiring that the DP-polynomials are primitive of each other  (in the sense of Remark \ref{integral}), i.e. we only require  that they satisfy  condition (1) of Definition \ref{G-modules}.

\begin{proposition}
\label{descente}
Let   $\mathcal H= \{H_1,\dots,H_{t+1}\}$ be a    finite set   of DP-polynomials satisfying condition (1) of Definition \ref{G-modules} with respect to a linear form $z $ of $R.$  Then $z^i \in \ann_R( H_i) $ for all $i=1, \dots, t+1. $

\noindent Let $I$ be an ideal of $R$ such  that  $\ann_R(H_{t+1})=I+(z^{t+1}). $ Assume  $z$  regular modulo $I$, then
$$\ann_R (H_{i})=I+(z^{i})$$ for all $i=1,\dots, t+1$. In particular $R/I$ is a one-dimensional Gorenstein $\res$-algebra and $I^{\perp} $ contains  $\{H_1,\dots,H_{t+1}\}.$
\end{proposition}
\begin{proof} First of all,  $z^i \in \ann_R H_i $ for all $i=1, \dots, t+1 $ since $z^i \circ H_i=0 $ by assumption.
We proceed now by descending recurrence.
Given $a\in I+(z^i)$, $i\le t$, we have $a z \in I+(z^{i+1})=\ann_R(H_{i+1})$.
Hence
$$
(a z)\circ H_{i+1}=a \circ (z \circ H_{i+1})= a \circ H_i=0
$$
so $ a\in \ann_R(H_i)$.
For all $ a\in \ann_R(H_i)$ we have $0= a \circ H_i= (a z) \circ H_{i+1}$.
Hence $a z \in  \ann_R(H_{i+1})= I + (z^{i+1})$.
Since $z$ is a non-zero divisor modulo $I$  we get $a\in I+(z^i)$. It follows that $R/I$ is one-dimensional Gorenstein $\res$-algebra because $I +(z) = \ann_R H_1$ is Artinian Gorenstein and   $z$ is  $I$-regular. We conclude by Theorem \ref{bijecGor}.
\end{proof}

\medskip

 We give an answer to Questions \ref{2} and \ref{1}  and we describe all the Gorenstein ideals $I$ of $R$ of multiplicity $e=e(R/I) $ such that $I^{\perp} $ contains a $G$-admissible set $\mathcal H= \{H_1,\cdots,H_{t+1}\}  $ with $t\ge e=\dim_{\res}\langle H_1\rangle. $  In particular this allows us  to find a system of generators of the Gorenstein ideal $I$ by a finite effective procedure. We note  that a procedure has been implemented in Singular using the next  result.

\begin{theorem}
\label{lift}
Let $\mathcal H= \{H_1,\cdots,H_{t+1}\} $ be a $G$-admissible set with respect to a linear form $z$ with $t\ge e=\dim_{\res}\langle H_1\rangle. $  Then the following conditions are equivalent:
\vskip 1mm
\noindent
$(i)$  There exists       a $G$-admissible $R$-submodule of $\Gamma$    with respect to $z\in R, $ say $ M_{\mathcal H},  $ obtained by  completion of $\mathcal H. $
 In particular  if  $I=\ann_R(M_{\mathcal H} ),    $  then $A=R/I$ is one-dimensional and Gorenstein, $z $ is regular modulo $I, $  and $\ann(H_{t+1})=I+(z^{t+1}).$
\vskip 1mm
\noindent
$(ii)$  $ \ann_R(H_{t+1}) = (h_1, \dots, h_r) +(z^{t+1}) $  where   $h_1, \dots, h_r$ are the elements  of a minimal  standard  basis    of $\ann_R(H_{t+1})$  with $v({h_i})  \le t  $ and   there exist  $\alpha_1,\cdots, \alpha_r\in R$
such that $z$ is regular modulo $(h_1+ \alpha_1 z^{t+1}, \dots, h_r+ \alpha_r z^{t+1}).  $
\vskip 1mm
\noindent
If this is the case,  then  $I=(h_1+ \alpha_1 z^{t+1}, \dots, h_r+ \alpha_r z^{t+1}) $ and $I^{\perp}=  M_{\mathcal H}.  $

 \end{theorem}
\begin{proof}
Assume $(i), $ then there exists       a $G$-admissible $R$-submodule of $\Gamma$    with respect to $z\in R, $ say $ M_{\mathcal H},  $ obtained by  completion of $\mathcal H.$
 Let  $I=\ann_R(M_{\mathcal H} ) $ be the corresponding ideal. By Theorem \ref{bijecGor},  $R/I$ is one-dimensional and Gorenstein,  $z $ is regular  modulo  $I, $ and
$\ann_R(H_{t+1})=I+(z^{t+1}). $
 By  \thmref{minimalGB}, given  $h_1, \dots, h_r$  the elements  of a minimal  standard  basis    of $\ann_R(H_{t+1})$  with $v({h_i})  \le t,   $ we know that there is a family of elements $\alpha_1,\cdots, \alpha_r\in R$
such that $\{h_1+ \alpha_1 z^{t+1}, \dots, h_r+ \alpha_r z^{t+1}\}$
 is a minimal standard  basis of $I$ and  hence  $z$ is a non-zero divisor modulo  $(h_1+ \alpha_1 z^{t+1}, \dots, h_r+ \alpha_r z^{t+1})=I. $
Conversely, if we  assume $(ii) $ and consider the ideal $I= (h_1+ \alpha_1 z^{t+1}, \dots, h_r+ \alpha_r z^{t+1}),  $    then    $\ann_R(H_{t+1})=(h_1,\dots, h_r) +(z^{t+1}) = I +(z^{t+1}).  $  Because  $z $ is regular  modulo  $I, $ then  by Lemma \ref{descente},   $ R/I$ is  Gorenstein one-dimensional,   hence  $I^{\perp} $ is $G$-admissible with respect to $z$ and it contains $\mathcal H. $ Hence $I^{\perp}=   M_{\mathcal H}.$
 \end{proof}

\begin{remark} \label{alpha}
We remark that the previous result is effective.
Given a finite set  $\mathcal H= \{H_1,\cdots,H_{t+1}\} $ with $t\ge e=\dim_{\res}\langle H_1\rangle$, by using \cite{E-InvSyst14} we can check if $\mathcal H$ is G-admissible  with respect to a linear form $z$.
We can   compute $ \ann_R(H_{t+1}) $ and determine      the elements   $h_1, \dots, h_r$  of a minimal  standard  basis    of $\ann_R(H_{t+1})$  with $v({h_i})  \le t$.
Let $\{h_1,\cdots,h_s\}$ be a minimal system of generators of  $J=(h_1, \dots, h_r).  $
Let $\alpha_1, \dots, \alpha_s$ be elements of $R  $ such that $z$ is a regular superficial element modulo  $I_{\alpha}= (h_1+ \alpha_1 z^{t+1}, \dots, h_s+ \alpha_s z^{t+1})$, then $R/I_{\alpha} $ and $R/J$ are Gorenstein one-dimensional rings with $I_{\alpha}^*=J^*$ sharing the same Hilbert function.
In fact $R/I_{\alpha} $ is Gorenstein one dimensional by Proposition \ref{descente} and we conclude concerning  $R/J$ by  Theorem \ref{minimalGB}.
  \end{remark}

\begin{example} \label{2ideals}
Consider $R=\res[\![x,y,z]\!].$ According to \remref{alpha},  we wish to construct  the   ideals $I$   in $R$ such that $R/I $ are  Gorenstein local rings  of dimension $1 $ and multiplicity $5 $ sharing the same $G$-admissible set with respect to $x\in R.$  Consider a polynomial $H_1 \in \Gamma=\res_{PD}[Y,Z] $ such that $\dim_k \langle H_1\rangle=5. $ Let $H_1= Z^{2} +Y^{3}. $
 One can   verify that
 \begin{multline*}
\mathcal H=\{ H_1, H_2=XH_1,  H_3=X^{2}H_1, H_4= X^{3}H_1+ Y^{4} Z +Y Z^{3},H_5=X H_4, H_6=X H_5\}
 \end{multline*}
is a finite G-admissible set with respect to $x.$   

The elements of a  minimal standard  basis of $  \ann_R(H_6)$ of valuation $\le e=5$ are $h_1=yz-x^3, h_2=z^2-y^3, h_3=y^4-x^3z.$ 
We verify that  $x$ is regular modulo the ideal $J=(h_1, h_2, h_3) = (h_1, h_2), $ hence by Theorem \ref{lift}  $R/J$ is one-dimensional Gorenstein.
Notice that the ideal $J=(h_1,h_2,h_3)=(h_1,h_2)$ is  the defining ideal in $R$ of the semigroup ring $\res[\![t^5,t^6,t^9]\!].$
Consider now  the ideal  $I_{\alpha}=(h_1 + \alpha_1x^6 ,h_2+\alpha_2x^6)$ with $\alpha_1, \alpha_2 \in R.$
Since $x$ is regular modulo $I_{\alpha}$ for every $\alpha_1, \alpha_2 \in R,  $ then   $I_{\alpha} $ describes all the ideals of $R$ such that $ R/I_{\alpha} $ is a one-dimensional Gorenstein ring of multiplicity
$5 $  and the dual is a completion of  $\mathcal H.$ In particular the family of Gorenstein ideals $I_{\alpha} $ share the same associated graded ring, hence the same Hilbert function,  because they have the same tangent cone.
\end{example}

In general, it is hard to prove that two ideals  are non analytically isomorphic.
In the next example we give two deformations $I_1$, $I_2$ of $J$, see \remref{alpha}, such that $I_1, I_2, J$ are pairwise non-analytically isomorphic and they come from the same $G$-admissible set.

\begin{example}
\label{non-iso}
Consider $R=\res[\![x,y,z]\!] $ and  the $G$-admissible set with respect to $x\in R$
$$
\mathcal H=\{ H_1= Y^3+Z^2, H_i=X^{i-1}H_1, i=2, \cdots ,6\}.
$$
Hence $H_6=X^5Y^3+X^5Z^2$ and $\ann_R(H_6)  = (yz-x^7, z^2-y^3, x^6).$
Following the notations of \remref{alpha}, $J=(h_1,h_2)$ with $h_1=yz-x^7, h_2=z^2-y^3$.
Notice that $J$ is the defining ideal of the Gorenstein monomial curve $\res[\![t^5,t^{14},t^{21}]\!].$
Consider the following two deformations of $J$:
$I_1=(h_1+x^6,h_2)$ and $I_2=(h_1,h_2+x^6)$,
that is respectively $(\alpha_1,\alpha_2)=(x^6,0)$ and  $(\alpha_1,\alpha_2)=(0,x^6).$
Both ideals define reduced  one-dimensional Gorenstein rings.
By using  Singular we compute the  multiplicity sequences of the above ideals: $\{5,5,4,1,...\}$ for $J$, $\{5,5,2,2,1,...\}$ for $I_1$, and $I_2$ defines an ordinary  singularity with $5$ non-singular branches.
Hence $J, I_1, I_2$ are pairwise non-analytically isomorphic.
\end{example}

\section{$0$-dimensional Gorenstein schemes}
 As an application of the previous results,  in this section we will present  some examples and discussions that we hope will be useful in the theory of Gorenstein linkage (G-linkage).  Note that the lack of a general structure of homogeneous  Gorenstein ideals of higher codimension   is the
main obstacle to extending the Gorenstein liaison theory in codimension at least three;
the codimension two Gorenstein liaison case is well understood,  see \cite{KMMNP01}.
See, for instance,
\cite{MP97}, \cite{KM83} and \cite{IS05} for some constructions of particular families of  Gorenstein algebras.

We say that two schemes $X$ and $Y$ in ${\mathbb P}^n_{\res}$  that are reduced
and without common components are directly Gorenstein-linked if their union is
arithmetically Gorenstein.
More   generally,   if we call $I(G) $ the Gorenstein ideal in $R=\res[x_0, \dots, x_n]$,    then
$$ I(G) : I(X) =I(Y)\ \ \  \text{and} \ \ \  I(G) : I(Y) =I(X). $$
Equivalence classes in Gorenstein linkage are determined by the equivalence relation generated by direct Gorenstein  linkage. 
 In terms of the inverse system,  it is easy to verify that $ I(G) : I(X) =I(Y)$ is equivalent to
$$I(X)^{\perp}  =   I(Y) \circ I(G)^{\perp}. $$

We say
$X $ is  {\it{glicci}} if it is in the Gorenstein linkage equivalence class of a complete intersection. In analogy  to the codimension 2 theory (where Gorenstein is complete intersection),
one hopes that every arithmetically Cohen-Macaulay scheme is glicci  and, in
particular, every finite set of points in $\bf {P}^3$  should be glicci, see for instance \cite{Har02a}, \cite{KMMNP01}.  This was verified by R. Hartshorne in \cite{Har02a} and by D. Eisenbud, R. Hartshorne and F.O. Schreyer in \cite{EHS15}  if $d=|X| \le 33, \ d=37, 38.$

Usually  it is very difficult to construct a set  of distinct points $G$ whose defining ideal is Gorenstein. General points are very far in general from being Gorenstein. The problem becomes  even more difficult if the set $G$ must contain  a given set $X$ of points as the linkage theory requires.  We take advantage of the results of this paper   for giving  explicit examples of reduced one dimensional Gorenstein $\res$-algebras $A=\res[x_0, \dots,x_n]/I(G)$ with  given Hilbert function. Our approach aims to be constructive, but for the moment we cannot deduce a theoretic  approach.

Our goal is the following. Given the integers $d; e$  and an appropriate Hilbert function of degree d+e admissible for a  Gorenstein graded $\res$-algebra of dimension one, we will construct some  examples of   reduced Gorenstein schemes $G$   in $\mathbb P^3_{\res}$ with the given Hilbert function.
We verify  that $G$  contains a degree $d$ subscheme $X$ with generic Hilbert function whose complement $Y$ of degree $e$ also has generic Hilbert function. It turns out that there is only a finite number of possibilities for the Hilbert function of such Gorenstein schemes. We recall that in our construction the Hilbert function  is determined by the suitable choice of $H_1.$  In fact the entries of the $h$-vector are determined by $\dim_k \langle H_1\rangle _j.$
 We say that a $0$-dimensional scheme in $\mathbb P^n_{\res}$ has generic Hilbert function if it is maximal, that is $\HF_X(j) = {\rm Max} \{ \deg X, \binom{n+j}{j}  \}.$ In general $I(X) $ is level, hence also the structure of $ I(X)^{\perp}$ is known by \cite{MT18}.
  Different methods for constructing points which are glicci in special cases are known.   Our hope is that our computations could suggest    some theoretical methods to use more abstractly.

In this section let   $R=\res[x,y,z,w] $   and hence $\Gamma=\res_{DP}[X,Y,Z,W]. $  We assume the ground field $\res$ is of characteristic zero.
 We may assume without loss of  generality that the points we want to construct do not lie on $w=0.$ Hence we want to find a  $G$-admissible set  $ \mathcal H=\{H_1,\dots, H_{t+1}\} $ in $\Gamma$ with respect to $w $ such that $\dim_k \langle H_1\rangle= d+e. $
Moreover if $t= $ maximum degree of the generators of $I(G)  $ or $t= \deg H_1, $ then by the previous results:
 $$I(G) = \ann_R(H_{t+1})_{\le t}R.$$
 If $X$ is   glicci, then there exists  a $G$-admissible sequence  $\mathcal H= \{H_1,\dots, H_{t+1}\} $ such that  $$(I(X) +(w^i))^{\perp} \subseteq \langle H_{i}\rangle  $$ for all $i=1, \dots, t+1. $


\begin{example}
\label{4points} ($5$ Gorenstein points)
We know that a set of {\bf{$4$ points in $\mathbb P^3_{\res}$ is glicci,}} in fact they are linked to a complete intersection by a Gorenstein scheme of length $5. $  This is a very special case, in fact  if we consider a set $X$ of $4$ points in linear general position and we take a sufficiently general point $Y$, then $G=X \cap Y$ consists of  $5$ points which are Gorenstein (but not a complete  intersection).
If we consider $H_1= X^2+Y^2+Z^2, $ then $\dim_k\langle H_1\rangle= \dim_k\langle 1, X,Y,Z, H_1\rangle=5  $ and we can take the $G$-admissible set $\mathcal H = \{ H_1, H_2= WH_1+ X^3+Y^3+Z^3, H_3= WH_2+ X^4+Y^4+Z^4\} $ with respect to $w \in R. $ In this case, 
then $I= \ann_R (H_3)_{\le 2}R $ is radical and it is the defining ideal of a Gorenstein set of $5$   points containing $4$ points in linear general position.
\end{example}

\begin{example}
\label{14points} (Gorenstein set of $14$ points)
It is known that  a set of  {\bf{$10$ points in $\mathbb P^3_{\res}$ is  glicci, }} in fact they are linked to a set of $4$ points by a Gorenstein scheme of length $14. $   In this example we use the inverse system to find a  Gorenstein zero-scheme consisting of  $d+e=14$ points which are the union of two reduced subschemes $X$ and $Y$ respectively of $d=10 $ and $e=4$ points with maximal Hilbert function.

Consider $H_1=X^3Y+Y^3Z+XZ^3 \in \res_{DP}[X,Y,Z] $ and we observe that  $\dim_k \langle H_1 \rangle=14 $ with Hilbert function $h_i= \dim_k \langle H_1\rangle_i:\ \ $  $ \{1, 3, 6, 3,1\}$.
 Recall that $\ann_R(H_1) = I(G) +(w). $ We want to construct a $1$-dimensional Gorenstein graded algebra $R/I(G) $ from a   $G$-admissible set  $\mathcal H =\{H_1, H_2,H_3, H_4\}. $ Notice that in Theorem \ref{lift}, we can choose $t=3$ because by the Hilbert function we know that the ideal is generated in degree $3.$ By using Singular's library \cite{E-InvSyst14},  given $H_j $ with $ j=1, \dots, t,$ it is possible to find the admissible elements $C_j \in \res_{DP}[X,Y,Z] $  such that $H_{j+1} =WH_j +C_j.$ In the first step   we prove that we can choose $C_1  $  any homogeneous form of degree $5$ in the variables $X,Y,Z.$ Take $C_1=X Y^2 Z^2.$

Consider
$H_2=WH_1+X Y^2 Z^2.$  {\it We realize computationally that this choice of $C_1$ fixes the successive elements}  $C_2= -2X^4Y^2-XY^4Z+2Z^6$ and $C_3= -XY^6-2X^5YZ-2X^2Y^3Z^2-2X^3Z^4, $
then  $H_4= W^3H_1+W^2 C_1+WC_2 + C_3.$ Hence  $$I= \ann_R(H_4)_{\le 3} R=(x^3-y^2z+z^2w, x^2y-z^3+2xw^2, $$ $$xy^2-xzw+2yw^2,   y^3-xz^2+yzw, x^2z+2w^3, xyz-y^2w+zw^2,yz^2-x^2w) $$   is   Gorenstein and $\dim R/I=1.$  Using  Singular,  we can verify that the ideal is radical,  hence it is the defining ideal of a Gorenstein set of $14$  distinct points. We can also prove  that it is the union of two subsets of points comprised  of $10$ and $4$ points with maximal Hilbert function.

It would be more interesting to have a theoretic argument proving that starting from a given subscheme $X$ of $d=10$ distinct points with maximal Hilbert function, then it is always possible to find $e=4$ distinct points, say $Y, $ such that the union defines a $0$-dimensional  Gorenstein scheme of length $d+e=14. $

   In our case $ I(Y) $ is generated by $6$ quadrics $q_1=x^2+wL_1, q_2=xy+wL_2, \dots, q_6=z^2+wL_6$ with $L_i$ suitable linear forms  in $R  $ and $I(X) = \ann_R(q_1 \circ H_4, \dots, q_6 \circ H_4)_{\le 3}R $ where $\mathcal H = \{H_1, H_2, H_3, H_4\} $   is a $G$-admissible set.
Our approach suggests to find $J= (q_1, \dots, q_6)=I(4)  $ imposing the following conditions:
$ (I(X) +(w^i))^{\perp} = J \circ H_i$, $i=1,\cdots, 4$, with $\mathcal H = \{H_1, H_2, H_3, H_4\} $   $G$-admissible.

In our particular case the system has a solution and  we have:
$L_1= 3 x +12 y + 3 z+ 7 w, L_2= 6 x + 11 y + 9 z+ 11 w, L_3= 12 x +10 y + 11 z+ 7 w, L_4= 9 x + 7 y + 5 z+ 11 w, L_5= x + 10y +4 z+ 3 w, L_6= 10 x + 2 y + 10 z+ 2 w $ and   $I(X)=(-x^3+3x^2y+3y^3-y^2z+3xz^2+3z^3+xyw,
-2x^3+6x^2y-4y^3-2y^2z-4xz^2+yz^2+6z^3+x^2w,
-3x^3-x^2y-2y^3-2y^2z-2xz^2-z^3+z^2w,
-6x^3-4x^2y+5y^3-6y^2z+6xz^2-4z^3+yzw,
-3x^3+x^2y+4y^3+xyz-3y^2z+4xz^2+z^3+y^2w
2x^3-3x^2y+xy^2-3y^3+2y^2z-3xz^2-3z^3+xzw).$
 \end{example}

 \begin{example}
\label{20points} (Gorenstein sets of $30$ and $55$  points)
This case was considered by Hartshorne in \cite{Har02a} to be a possible counterexample to the ``conjecture". But Eisenbud, Hartshorne and Schreyer in  \cite{EHS15} proved that  a set of  {\bf{$20$ points in $\mathbb P^3_{\res}$ is  glicci, }} in fact they are linked to a set of $10$ points by a Gorenstein scheme of length $30. $   In this example we use the inverse system to find a set $G$ of $d+e=30$ Gorenstein points which are the union of two reduced subschemes $X$ and $Y$ respectively of $d=20 $ and $e=10$ points with maximal Hilbert function.  We may proceed as in the previous example.

In this case $\deg H_1= 6$ because the socle degree of $R/I(G) $ is $6.$
 We consider $H_1=X^6 + Y^6+Z^6 + X^5Y+Y^5Z+XZ^5 \in \res_{DP}[X,Y,Z] $ because $\dim_k \langle H_1\rangle =30 $ and the Hilbert function is: $\{1,  3, 6, 10,  6, 3, 1\}$.
 Recall that $\ann_R(H_1) = I(G) +(w). $ We want to construct a $1$-dimensional Gorenstein graded algebra $R/I(G) $ from a    $G$-admissible set  $\mathcal H =\{H_1, \dots, H_5\}. $ Notice that in Theorem \ref{lift}, we can choose $t=4$ because by the Hilbert function we know that the ideal is generated in degree $4.$ By using Singular's library \cite{E-InvSyst14}  we can find $C_1, \dots, C_4$  and we get $H_5=W^4 H_1+ W^3 C_1+ \dots +WC_3+C_4.$ At the first step we choose $C_1= X^2Y^2Z^2, $ then as before we realized that  all  the successive forms are uniquely determined.

 In the same paper, Eisenbud, Hartshorne and Schreyer stated that the degree of the smallest collection of general points in
 $\mathbb P^3_{\res}  $  not yet known to be glicci is $34.$  As before, what we are able   to construct, with our methods, is  a set $G$ of $55$ Gorenstein points containing a subscheme of $34$ distinct points with generic Hilbert function,  whose complementary   $21$ points has also generic Hilbert function.

In this case $\deg H_1= 8$ because the socle degree of $R/I(G) $ is $8.$ We consider $H_1=X^8 + Y^8+Z^8 + X^3Y^3Z^2+X^2Y^3Z^3+X^3Y^2Z^3 \in \res_{DP}[X,Y,Z] $ because $\dim_k \langle H_1\rangle =55 $ and its Hilbert function is: $\{ 1, 3, 6, 10, 15, 10,  6, 3, 1\}.$
 Recall that $\ann_R(H_1) = I(G) +(w). $ We want to construct a $1$-dimensional Gorenstein graded algebra $R/I(G) $ from a    $G$-admissible set  $\mathcal H =\{H_1, \dots, H_6\}. $ Notice that in Theorem \ref{lift}, we can choose $t=5$ because by the Hilbert function we know that the ideal is generated in degree $\le 5.$ By using Singular's library \cite{E-InvSyst14} we can find $C_1, \dots, C_5$  and we get $H_6=W^5 H_1+ W^4 C_1+ \dots +WC_4+C_5.$ We choose $C_1= XY^4Z^4, $ then by the algorithm we realize that, as in the previous examples,  all the successive ``constants" $C_i$  are uniquely determined. For the moment we cannot understand the reason for  this constraint. We think that this could be the key of a constructive general argument.
\end{example}

\section{Gorenstein Semigroup rings}

Semigroup rings are a broad class of one-dimensional domains and they are a test case for many open problems. Hence we end this paper with the construction of the inverse system of the defining  ideal of a family of semigroup rings. In particular we study a corresponding $G$-admissible set   in  $\Gamma$ according  to Theorem \ref{lift}.  Let $b \ge 2 $ be an integer. Consider  
$$A= \res[[ t^{3b}, t^{3b+1}, t^{6b+3}]].$$
It is easy to see that   $A= \res[\![x,y,z]\!]/I $ where $I= (xz-y^3, z^b-x^{2b+1}). $
 Thus $A$ is a one-dimensional Gorenstein local domain of type $(2,b)$ and multiplicity $e=3b,$ see  \cite{GHK07}, Example 5.5.
It is an interesting class of ideals because  the associated graded ring is not Cohen-Macaulay even if the ideal is a complete intersection.  The Hilbert function of $A$   has a very special shape.
For every $b \ge 2 $ the Hilbert function of $A$ has    $b-1$ flats \cite{ERV14},
\begin{equation} \label{b-1flats} \HF_A(t)=
\begin{cases}
     1 & \ \ \text{$t=0$}, \\
    2 t+2 & \ \ \text{$t=1,\cdots, b-1$},\\
    2 b & \ \ \text{$t=b$},\\
    2b+1 & \ \ \text{$t=b+1$},\\
    2b+k  & \ \ \text{$t=b+2k, \ \ k=1,\cdots, b-1$},\\
    2b+k+1  & \ \ \text{$t=b+2k+1, \ \ k=1,\cdots, b-1$},\\
    3b  & \ \ \text{$t\ge 3b-1 $}.\\
\end{cases}
\end{equation}

\noindent
We  compute the inverse system $I^\perp=\langle H_t, t\in \mathbb N_{+}\rangle$ of $I$ with respect to $x$. Notice that the ideal $I$  is homogeneous in the weighted space defined by $\deg_b=(3b, 3b+1, 6b+3)$, i.e.
$$
\deg_b(x^iy^jz^k)=3bi+(3b+1)j+(6b+3)k.
$$
Because the ideal $I$ is homogeneous with respect to $\deg_b$,   the DP-polynomial $ H_t$ is homogeneous as well.  Notice that $e=3b=\dim_{\res}\langle H_1\rangle$.
We determine  that $\mathcal H= \{H_1,\cdots,H_{3b+1}\} $ is a G-admissible set in $\Gamma= \res_{DP} [X, Y, Z] $ with respect to $x  $ which determines $I.$

\medskip
\noindent {\bf{Claim}}:
$H_1=Y^2Z^{b-1}$ and $\deg_b(H_t)=6 b^2+ 3bt-1$, $t\ge 1$.

\noindent
In fact we know that $H_1$ is a generator of $(I+(x))^{\perp}$.
Since $I+(x)=(x,y^3, z^b)$, an easy computation give us $H_1=y^2z^{b-1}$.
Recall that $H_{t}=x\circ H_{t+1}. $ From this we get the second part of the claim.

We write $\delta(b,t)=6 b^2+ 3bt-1$, for $b\ge 2$ and $t\ge 1$;
we denote by $\mathcal M_{b,t}$  the set of homogeneous monomials of degree $\delta(b,t)$ in $\Gamma= \res_{DP} [X, Y, Z] $.
For all $t\ge 1$ consider the DP-polynomial of degree  $\delta(b,t)$
$$
H_t=\sum_{\it m\in \mathcal M_{b,t}} \it m.
$$

\medskip
\noindent {\bf{Claim}}:  $\ann_R(H_{3b+1})= I + (x^{3b+1}) $.

\medskip
\noindent
Notice that we only have to prove that
$$
\left\{
  \begin{array}{ll}
    (1) & x^{3b+1}  \circ H_{3b+1} =0 \\
    (2) & (xz-y^3)\circ  H_{3b+1}=0\\
    (3) & (z^b-x^{2b+1})\circ H_{3b+1}=0
  \end{array}
\right.
$$
\noindent
We set $\delta=\delta(b,3b+1)=15 b^2+3 b-1. $ Recall that $\deg(H_{3b+1})=\delta$.

\noindent
$(1)$
For all monomials $m=X^iY^jZ^k$  of degree $\delta$ we have to prove that $i< 3b+1$.
Assume that $i\ge 3b+1$.
The Diophantine equation
$$
(3b+1)j+(6b+3)k= \delta-3bi
$$
has the following solutions:
$$
(j,k)=(-2,1)(\delta-3bi) + \lambda (6b+3, -(3b+1))
$$
for all $\lambda \in \mathbb Z$.
Since $j, k\ge 0$ we deduce $\delta-3bi\ge 0$ and
$$
A=\frac{2(\delta-3bi)}{6b+3} \le \lambda \le B=\frac{\delta-3bi}{3b+1}.
$$
In particular $i< 5b+1$ and then $3b+1\le i< 5b+1$.
We can perform the corresponding  Euclidean divisions
$$
A= 5b-i+r_1, B=5b-i+r_2
$$
with $0< r_1=\frac{3i-9b-2}{6b+3}< r_2=\frac{i-2b-1}{3b+1}<1$.
Hence there are no integers $\lambda$ such that $A\le \lambda \le B$.
Hence $i < 3b+1$.

\noindent
(2) We  prove that $xz\circ \mathcal M_{b,3b+1}= y^3\circ \mathcal M_{b,3b+1}$.
Given  $\it m=X^iY^jZ^k\in \mathcal M_{b,3b+1}$  we have that $xz\circ {\it m} =0$ if $i=0$ or $k=0$.
Hence
$xz\circ \mathcal M_{b,3b+1}$ is the set of monomials ${\it m'}= X^{i-1} Y^j Z^{k-1}$ such that ${\it m} =X^iY^jZ^k\in \mathcal M_{b,3b+1}$ and $i,k \ge 1$.
Since ${\it m'}=y^3 \circ (Y^3 {\it m'})$ and  $y^3 {\it m'}\in \mathcal M_{b,3b+1}$ we get
$xz\circ \mathcal M_{b,3b+1}\subset  y^3\circ \mathcal M_{b,3b+1}$.
Given  $\it m=X^iY^jZ^k\in \mathcal M_{b,3b+1}$  we have that $y^3\circ {\it m} =0$ if $j \le 2$.
Hence
$y^3\circ \mathcal M_{b,3b+1}$ is the set of monomials ${\it m'}= X^{i} Y^{j-3} Z^{k}$ such that ${\it m} =X^iY^jZ^k\in \mathcal M_{b,3b+1}$ and $j\ge 3$.
Since ${\it m'}=xz \circ (XZ {\it m'})$ and  $XZ {\it m'}\in \mathcal M_{b,3b+1}$ we get
$y^3\circ \mathcal M_{b,3b+1}\subset  xz\circ \mathcal M_{b,3b+1}$.

\noindent
(3) We   prove that $z^b\circ \mathcal M_{b,3b+1}= x^{2b+1}\circ \mathcal M_{b,3b+1}$.
We proceed as for $(2)$.
\vskip 2mm
Similarly one can prove that $\ann_R(H_{t}) = I +(x^{t}) $ for every $t \ge 1.$

\bigskip
\noindent
For instance if  $b=3$ we get that
$H_{10}=x^2y^2z^5+xy^5z^4+y^8z^3+x^9y^2z^2+x^8y^5z+x^7y^8. $
In this case  $$\ann_R(H_{10})=   (xz-y^3, z^3-x^{7}) + (x^{10}) = I + (x^{10}). $$


\end{document}